\documentclass[12pt]{article}

\usepackage[cp1251]{inputenc}
\usepackage[T2A]{fontenc}
\usepackage[english, ukrainian]{babel}
\usepackage{amsmath,amsfonts,amssymb,amsthm,cite}
\usepackage{geometry}
\theoremstyle{plain}
\newtheorem{theorem}{Theorem}

\newtheorem{corollary}{Corollary}
\newtheorem{definition}{Definition}
\theoremstyle{definition}

\newcommand{\nonprint}[1]{}

\begin{document}
\selectlanguage{english}

\begin{flushleft}

\Large
\begin{center}
\textbf{Differential systems in Sobolev spaces with generic inhomogeneous boundary conditions}
\end{center}

\begin{center}
\textbf{Vladimir Mikhailets, Olena Atlasiuk}
\end{center}

\Large


\end{flushleft}

\normalsize

\begin{abstract}
The paper contains a review of results on linear systems of ordinary differential equations of an arbitrary order on a finite interval with the most general inhomogeneous boundary conditions in Sobolev spaces. The character of the solvability of such problems is investigated, their Fredholm properties are established, and their indexes and the dimensions of their kernels and co-kernels are found. In addition, necessary and sufficient conditions of continuity in the parameter of the solutions of the introduced classes of boundary-value problems in Sobolev spaces of an arbitrary order are obtained.
\end{abstract}


\vspace*{+1em}

\verb"2020 Mathematics Subject Classification: 34B05, 34B08, 34B10, 47A53"

\vspace*{+1em}

\section{Introduction}
The study of systems of ordinary differential equations is an important part of many investigations in modern analysis and its applications (see, for example, \cite{BochSAM2004} and references therein). Unlike the Cauchy problem, the solutions of inhomogeneous boundary-value problems for differential equations/systems may not exist and/or may not be unique. Therefore, the question of the character of the solvability of such problems is fundamental for the theory of differential equations. Thus, Kiguradze \cite{Kigyradze1987, Kigyradze1975} and Ashordia \cite{Ashordia} investigated systems of first-order differential equations with general inhomogeneous boundary conditions of the form
\begin{equation}\label{kig}
y^{\prime}(t) + A(t)y(t)=f(t), \quad t\in(a,b), \quad By=c.
\end{equation}
Here, the $m \times m$ -- matrix-valued function $A(\cdot)$ is Lebesgue integrable over the finite interval $(a, b)$; the vector-valued function $f(\cdot)$ belongs to Lebesgue space $L\left((a, b); \mathbb{ R}^m\right)$; the vector $c \in \mathbb{R}^{m}$, and $B$ is an arbitrary linear continuous operator from the Banach space $C\left([a,b]; \mathbb {R}^{m}\right)$ to $\mathbb{R}^{m}$, with arbitrary $m \in \mathbb{N}$. The boundary condition in \eqref{kig} covers the main types of classical boundary conditions; namely: Cauchy problems, two-point and multipoint problems, integral and mixed problems. The Fredholm property with zero index was established for problems of the form \eqref{kig}. Moreover, the conditions for the problems to be well posed were obtained. The limit theorem for the sequence of their solutions in the space of continuous vector-valued functions on $[a, b]$ is proved.

These results were further developed in a series of articles by Mikhailets and his disciples \cite{GnypKodMik15, GKM2017, KodlyukMikR2013, KM2013, MikhailetsChekhanova, MMS2016,  MikPelRev2018}. Specifically, they studied the systems of differential equations of an arbitrary order $r \in \mathbb{N}$. In this case, the operator $B$ specifying the inhomogeneous boundary condition is an arbitrary linear continuous operator from the complex Banach space $C^{r-1}\left ([a,b]; \mathbb{C}^{m}\right)$ to $\mathbb{C}^{rm}$. They obtained conditions for the boundary-value problems to be well posed and proved limit theorems for the sequences of solutions of such problems in the space $C^{r-1}\left ([a,b]; \mathbb{C}^{m}\right )$. These results significantly generalize Kiguradze's theorems even in the $r=1$ case. Moreover, limit theorems for Green's matrices of such boundary-value problems were established for the first time \cite{KodlyukMikR2013, MikhailetsChekhanova}. These results have already found application to the analysis of multipoint boundary-value problems \cite{Atl3}, as well as to the spectral theory of differential operators with distributions in coefficients  \cite{GoriunovMikhailetsPankrashkin2013EJDE, GoriunovMikhailets2010MFAT, GoriunovMikhailets2012UMJ, GoriunovMikhailets2010MN}.

Note that boundary-value problems with inhomogeneous boundary conditions containing derivatives of the unknown vector-valued function of integer and$/$or fractional orders that can be equal to or greater than the order of the differential equation naturally arise in some mathematical models \cite{Kilbas2006, Kra1961, Luo1991, Vent1959}). The theory of such problems has not been developed yet and contains few results even for the case of ordinary differential equations. The study of such problems requires new approaches and methods of the analysis \cite{GKM2017, MMS2016}.

The aim of this article is to give a brief survey of this theory for linear systems of ordinary differential equations of an arbitrary order with the most general (generic) inhomogeneous boundary conditions in Sobolev spaces.

Let a finite interval $(a,b)\subset\mathbb{R}$ and parameters $$\{m, \, n+1, \,r,\,l\} \subset \mathbb{ N}, \, 1\leq p\leq \infty$$ be given. Let \begin{align*}
&W_p^{n+r}\bigl([a,b];\mathbb{C})\\
&:= \bigl\{y\in C^{n+r-1}([a,b];\mathbb{C})\colon y^{(n+r-1)}\in AC[a,b], \, y^{(n+r)}\in L_p[a,b]\bigr\}
\end{align*}
denote the corresponding complex Sobolev space; set $W_p^0:=L_p$. This space is Banach with respect to the norm
$$
\bigl\|y\bigr\|_{n+r,p}=\sum_{k=0}^{n+r}\bigl\|y^{(k)}\bigr\|_{p},
$$
with $\|\cdot\|_p$ standing for the norm in the Lebesgue space $L_p\bigl([a,b]; \mathbb{C}\bigr)$. Similarly, we let $(W_p^{n+r})^{m}:=W_p^{n+r}\bigl([a,b];\mathbb{C}^{m}\bigr)$ and $(W_p^{n+r})^{m\times m}:=W_p^{n+r}\bigl([a,b];\mathbb{C}^{m\times m}\bigr)$ denote the Sobolev spaces of vector-valued functions and matrix-valued functions, res., whose ele\-ments belong to the function space $W_p^{n+r}$.

We consider the following linear boundary-value problem:
\begin{equation}\label{bound_pr_1}
(Ly)(t):=y^{(r)}(t) + \sum\limits_{j=1}^rA_{r-j}(t)y^{(r-j)}(t)=f(t), \quad t\in(a,b),\\
\end{equation}
\begin{equation}\label{bound_pr_2}
By=c,
\end{equation}
where the matrix-valued functions $A_{r-j}(\cdot)\in (W_{p}^{n})^{m\times m}$, vector $c\in\mathbb{C}^{l}$, vector-valued function $f(\cdot)\in (W_{p}^{n})^{m}$,  and linear continuous operator
\begin{equation}\label{oper_B_class}
B\colon (W^{n+r}_p)^m\rightarrow\mathbb{C}^{l}
\end{equation}
are arbitrarily chosen; whereas the vector-valued function $y(\cdot)\in (W_{p}^{n+r})^m$ is unknown. If $l<rm$, then the boundary conditions are underdetermined. If $l>rm$, then the boundary conditions are overdetermined.

The boundary condition \eqref{bound_pr_2} consists of $l$ scalar conditions for the system of $m$ differential equations of the $r$ order. We represent vectors and vector-valued functions in the form of columns. The solution to the boundary-value problem \eqref{bound_pr_1}, \eqref{bound_pr_2} is understood as a vector-valued function $y(\cdot)\in (W_{p}^{n+r})^m$ satisfying equation \eqref{bound_pr_1} (for $n\geq 1$ everywhere, and for $n= 0$ almost everywhere) on $(a,b)$, and equality \eqref{bound_pr_2}. If the parameter $n$ increases, so does the class of linear operators \eqref{oper_B_class}. When $n=0$, this class contains all the operators that specify the general boundary conditions.

The solutions of equation \eqref{bound_pr_1} fill the space $(W_{p}^{n+r})^m$ if its right-hand side $f(\cdot)$ runs through the space $(W_{p}^{n})^m$. Therefore, boundary condition \eqref{bound_pr_2} with continuous operator \eqref{oper_B_class} is the most general condition for this equation.

It is known \cite{Ioffe} that, if $1\leq p < \infty$, then every operator \eqref{oper_B_class} admits the unique analytic representation
\begin{equation}\label{st anal}
By=\sum _{i=0}^{n+r-1} \alpha_{i}\,y^{(i)}(a)+\int_{a}^b \Phi(t)y^{(n+r)}(t){\rm d}t, \quad y(\cdot)\in (W_{p}^{n+r})^{m},
\end{equation}
for certain number matrices $\alpha_{s} \in \mathbb{C}^{rl\times m}$ and a  matrix-valued function $$\Phi(\cdot)\in L_{p'}([a, b]; \mathbb{C}^{rl\times m});$$ as usual, $1/p + 1/p'=1$. If $p=\infty$, this formula also defines a bounded operator $$B\colon (W_{\infty}^{n+r})^{m} \rightarrow \mathbb{C}^{rl}.$$ However, there exist other operators of this class generated by integrals over finitely additive measures. Hence, unlike  $p<\infty$ \cite{Atl2, GKM2017, KM2013},  the case of $p=\infty$ contains additional analytical difficulties.

The article is structured as follows.

Section \ref{section2} discusses the analysis of the solvability of the inhomogeneous boundary-value problem in the corresponding Sobolev spaces.

Section \ref{section3} gives examples that apply to the results of Section \ref{section2} and demonstrate the constructive character of these results.

Section \ref{section4} contains a limit theorem for the sequence of characteristic matrices of the considered boundary-value problems and some related results.

Section \ref{section5} contains definitions and necessary and sufficient conditions for the continuity of solutions to the boundary-value problems in a number parameter included in the coefficients of differential systems and boundary conditions.

Chapter \ref{section6} contains limit theorems for solutions of inhomogeneous multipoint boundary-value problems in separable and nonseparable Sobolev spaces.

\section{Solvability and characteristic matrix} \label{section2}

We rewrite the inhomogeneous boundary-value problem \eqref{bound_pr_1}, \eqref{bound_pr_2} in the form of a linear operator equation
\[
(L,B)y=(f,c).
\]
Here, $(L,B)$ is a bounded linear operator on the pair of Banach spaces
\begin{equation}\label{(L,B)}
(L,B)\colon (W^{n+r}_p)^m\rightarrow (W^{n}_p)^m\times\mathbb{C}^l,
\end{equation}
which follows from the definition of the Sobolev spaces involved and from the fact that $W_p^n$ is a Banach algebra.

Let $E_{1}$ and $E_{2}$ be Banach spaces. A linear bounded operator $T\colon E_{1}\rightarrow E_{2}$ is called a Fredholm one if its kernel and co-kernel are finite-dimen\-si\-onal. If $T$ is a Fredholm operator, then its range $T(E_{1})$ is closed in $E_{2}$, and its index
$$
\mathrm{ind}\,T:=\dim\ker T-\dim\big(E_{2}/T(E_{1})\big)\in \mathbb{Z}
$$
is finite (see, e.g., \cite[Lemma~19.1.1]{Hermander1985}).

\begin{theorem}\label{th_fredh high}
The bounded linear operator \eqref{(L,B)} is a Fredholm one with index $rm-l$.
\end{theorem}

The proof of Theorem \ref{th_fredh high} uses the well-known theorem on the stability of the index of a linear operator with respect to compact additive perturbations \cite{Atl1}.

Theorem \ref{th_fredh high} naturally raises the question of finding  $d$-characteristics of the operator $(L, B)$, i.e. $\operatorname{dim} \operatorname{ker}(L, B)$ and $\operatorname{dim} \operatorname{coker}(L, B)$. This is a quite difficult task because the Fredholm numbers may vary even with arbitrarily small one-dimensional additive perturbations.

To formulate the following result, let us introduce some notation and definitions.

For each number $i \in \{1,\dots, r\}$, we consider the family of matrix Cauchy problems:
\begin{equation}\label{zad kosh1}
Y_i^{(r)}(t)+\sum\limits_{j=1}^rA_{r-j}(t)Y_i^{(r-j)}(t)=O_{m},\quad t\in (a,b),
\end{equation}
with the initial conditions
\begin{equation}\label{zad kosh2}
Y_i^{(j-1)}(a) = \delta_{i,j}I_m,\quad j \in \{1,\dots, r\},
\end{equation}
where $Y_i(\cdot)$ is an unknown $(m\times m)$ -- matrix-valued function. As usual, $O_{m}$ stands for the zero $m\times m$ matrix, $I_{m}$ denotes the identity $(m\times m)$ -- matrix, and $\delta_{i,j}$ is the Kronecker delta. Each Cauchy problem \eqref{zad kosh1}, \eqref{zad kosh2} has a unique solution $Y_i\in(W_p^{n+r})^{m\times m}$ due to \cite[Lemma~4.1]{AtlMikh2024}). Certainly, if $r=1$, we use the designation  $Y(\cdot)$ for $Y_1(\cdot)$.

Let $\left[BY_i\right]$ denote the number $(l\times m)$ -- matrix whose $j$-th column is the result of the action of $B$ on the $j$-th column of the matrix-valued function~$Y_i$.

\begin{definition}\label{defin_harm}
A block rectangular number matrix
\begin{equation}\label{matrix_BY}
M(L,B):=\big(\left[BY_1\right],\dots,\left[BY_{r}\right]\big) \in \mathbb{C}^{l\times rm}
\end{equation}
is called the characteristic matrix of the inhomogeneous boundary-value problem \eqref{bound_pr_1}, \eqref{bound_pr_2}. This matrix consists of $r$ rectangular block columns $\left[BY_k\right]\in \mathbb{C}^{m\times l}$.
\end{definition}

Here, $mr$ is the number of scalar differential equations of the system \eqref{bound_pr_1}, and $l$ is the number of scalar boundary conditions in \eqref{bound_pr_2}.

\begin{theorem}\label{th dimker}
The dimensions of the kernel and co-kernel of the operator \eqref{(L,B)} are equal to the dimensions of the kernel and co-kernel of the characteristic matrix \eqref{matrix_BY}, respectively; i.e.,
\begin{gather*}\label{dimker}
\operatorname{dim}\operatorname{ker}(L,B)=
\operatorname{dim}\operatorname{ker}M(L,B),\\
\label{dimcoker1}
\operatorname{dim} \operatorname{coker}(L,B)=\operatorname{dim} \operatorname{coker}M(L,B).
\end{gather*}
\end{theorem}

\begin{corollary}\label{invertible}
The operator \eqref{(L,B)} is invertible if and only if $l=rm$ and the square matrix $M(L,B)$ is nonsingular.
\end{corollary}

In the $r=1$ case, Theorem~\ref{th_fredh high} and Corollary~\ref{invertible} are proved in~\cite{Atl1}. In the case where $l=rm$ and $p<\infty$, Corollary~\ref{invertible} is proved in~\cite{GnypKodMik15}. Theorem~\ref{th dimker} is also new for the systems of first order differential equations.

In Sobolev-Slobodetskii spaces, similar results for systems of first order differential equations were obtained in \cite{MikhSkor2021}.

The results concerning Theorems \ref{th_fredh high}, \ref{th dimker}, and Corollary \ref{invertible} given in this section were obtained in  \cite{AtlMikh2024}.

\section{Examples} \label{section3}

If all coefficients of the differential expression $L$ are constant, then the characteristic matrix of the corresponding boundary-value problem can be explicitly found in many instances (see, e.g., \cite{Gan1959}). In this case, the characteristic matrix is an analytic function of a certain square number matrix and coincides hence with some polynomial of this matrix.

\textbf{Example 1.} Consider the linear one-point boundary-value problem for first order constant-coefficient differential equation
\begin{equation}\label{1.6.1t1}
 (Ly)(t):= y'(t)+Ay(t)=f(t),\quad
t \in(a,b),
\end{equation}
\begin{equation}\label{1.3t1}
By= \sum _{k=0}^{n-1} \alpha_{k} y^{(k)}(a)=c,
\end{equation}
\noindent where $A$ is a constant $(m \times m)$ -- matrix; the vector-valued function $f(\cdot)$ belongs to the space~$(W_{p}^{n-1})^{m}$; the matrices $\alpha_{k}$ belong to the space $\mathbb{C}^{l\times m}$; $c \in \mathbb{C}^{l}$; the operators
\begin{equation*}\label{3.BY1}
B\colon (W_{p}^{n})^{m} \rightarrow\mathbb{C}^{l} \quad \mbox{and} \quad (L,B)\colon (W^{n}_p)^m\rightarrow (W^{n-1}_p)^m\times\mathbb{C}^l
\end{equation*}
act continuously, and $y(\cdot)\in (W_{p}^{n})^m$.

Let $Y(\cdot)\in (W_p^n)^{m\times m}$ denote a unique solution of the linear homogeneous matrix equation of the form \eqref{1.6.1t1} with the initial condition at the point $a$, namely:
\begin{equation*}\label{r31}
  Y'(t)+A Y(t)=O_{m},\quad t\in (a,b), \quad Y(a)=I_{m}.
  \end{equation*}

Put
\begin{equation*}\label{3.BY1}
M(L,B)=[BY]:=\left( B \begin{pmatrix}
                                              y_{1,1}(\cdot) \\
                                              \vdots \\
                                              y_{m,1}(\cdot) \\
                                            \end{pmatrix}
,\ldots,
                                    B \begin{pmatrix}
                                              y_{1,m}(\cdot) \\
                                              \vdots \\
                                              y_{m,m}(\cdot) \\
                                            \end{pmatrix}\right) \in \mathbb{C}^{m\times l}.
\end{equation*}
Then the fundamental matrix of system \eqref{1.6.1t1} and its $k$-th derivative have the following form:
\begin{gather*}
Y(t)= \operatorname{exp}\big(-A(t-a)\big), \quad Y(a) = I_{m}; \\
Y^{(k)}(t)= (-A)^k \operatorname{exp}\big(-A(t-a)\big), \quad Y^{(k)}(a) = (-A)^k, \quad k \in \mathbb{N}.
\end{gather*}

Substituting these value into the equality \eqref{1.3t1}, we have
$$M(L,B)=\sum_{k=0}^{n-1}\alpha_{k}(-A)^k.$$

Theorem \ref{th_fredh high} implies that $\operatorname{ind} \, (L, B)=\operatorname{ind} \, (M(L, B))= m-l$.

Therefore, by Theorem \ref{th dimker}, we obtain
\begin{gather*}
\operatorname{dim} \operatorname{ker}(L,B)=\operatorname{dim} \operatorname{ker}\left(\sum_{k=0}^{n-1}\alpha_{k}(-A)^k\right)=
m-\operatorname{rank}\left(\sum_{k=0}^{n-1}\alpha_{k}(-A)^k\right),  \\
\operatorname{dim} \operatorname{coker}(L,B)=-m+l+\operatorname{dim} \operatorname{coker}\left(\sum_{k=0}^{n-1}\alpha_{k}(-A)^k\right)=\\
l-\operatorname{rank}\left(\sum_{k=0}^{n-1}\alpha_{k}(-A)^k\right).\label{dimcoker}
\end{gather*}

From these formulas it follows that the $d$-characteristics of the problem do not depend on the length of the interval $(a, b)$.

\textbf{Example 2.} Let us consider a multipoint boundary-value problem for the system of differential equations \eqref{1.6.1t1}, with $A(t) \equiv O_{m}$. The boundary conditions at the points $\{t_k\}_{k=0}^N \subset [a, b]$ contain derivatives of integer and$/$or fractional orders (in the sense of Caputo \cite{Kilbas2006}). They have the next form
\begin{equation*}\label{3.BY1}
By=\sum _{k=0}^N \sum _{j=0}^s \alpha_{kj} y^{(\beta_{kj})}(t_{k})=c.
\end{equation*}
Here, the number matrices $\alpha_{kj} \in \mathbb{C}^{l \times m}$. The nonnegative numbers $\beta_{kj}$ are such that
\begin{equation*}\label{3.BY1}
\beta_{k,0}=0 \quad \mbox{for all} \quad k \in \{1, 2, \ldots, N\}. 
\end{equation*}

Theorem \ref{th_fredh high} implies that the index of the operator $(L, B)$ is equal to $m-l$.

Let us find the dimensions of its kernel and co-kernel. In this case,  $Y(\cdot)=I_{m}$. Therefore, the characteristic matrix has the form
\begin{equation*}\label{3.BY1}
M(L, B)=[BY]=\sum _{k=0}^N \sum _{j=0}^s \alpha_{kj} I_{m}^{(\beta_{kj})}=\sum _{k=0}^N \alpha_{k,0},
\end{equation*}
since the derivatives $\bigl(^{\textit{C}}\textmd{D} _{a+}^{\beta_{k,j}}I_{m}\bigr)=0$ if $\beta_{kj} >0$. Hence, according to Theorem \ref{th dimker}, we get
\begin{gather*}
\operatorname{dim} \operatorname{ker}(L,B)=\operatorname{dim} \operatorname{ker}\left(\sum_{k=0}^{N}\alpha_{k,0}\right)=m -
\operatorname{rank}\left(\sum_{k=0}^{N}\alpha_{k,0}\right),  \\
\operatorname{dim} \operatorname{coker}(L,B)= - m+l + \operatorname{dim} \operatorname{coker}\left(\sum_{k=0}^{N}\alpha_{k,0} \right)=l -
\operatorname{rank}\left(\sum_{k=0}^{N}\alpha_{k,0} \right).
\end{gather*}

It follows from these formulas that the $d$-characteristics of the problem do not depend on the choice of the interval $(a, b)$, points $\{t_k\}_{k=0}^N \subset [a, b]$, and matrices $\alpha_{kj}$, with $j \geq 1$.

\textbf{Example 3.} Consider a two-point boundary-value problem for a system of second-order differential equations generated by the expression
\begin{equation*}
    Ly(t):= y^{\prime \prime} (t)+Ay^{\prime}(t),\quad
t \in(a,b),
\end{equation*}
where $A$ is a constant matrix, with the boundary operator
\begin{equation*}
By=\sum^{n+1} _{k=0} \left(\alpha_{k} y^{(k)}(a)+ \beta_{k}y^{(k)}(b)\right).
\end{equation*}
Here, $\alpha_{k}$, $\beta_{k}$ are some rectangular number matrices. Then we have the continuous operator
$$
(L,B)\colon (W^{n+2}_p)^m\rightarrow (W^{n}_p)^m\times\mathbb{C}^l
$$
and characteristic matrix
$M(L,B) \in \mathbb{C}^{2m\times l}$.

It is easy to verify in this case that
$$ Y_1(t)\equiv I_m, \quad Y_2(t)=\varphi (A, t),$$
where, for each fixed $t \in [a, b]$, the function $\varphi (\lambda, t):=1-\exp(-\lambda(t-a))\lambda^{-1}$ is an entire analytic function of the variable $\lambda \in \mathbb{C}$.

Then
$$ [BY_1]=\sum_{k=0}^{n+1} \left(\alpha_kI_m^{(k)}(a)+\beta_kI_m^{(k)}(b)\right)=(\alpha_0+\beta_0)I_m,$$
$$ [BY_2]=\sum_{k=0}^{n+1} \left(\alpha_k\varphi^{(k)}(A, a)+\beta_k\varphi^{(k)}(A, b)\right).$$

But
$$Y_2^{(k)}(t)=(-1)^kA^k\exp (-A(t-a)), \quad k \in \{0, \ldots, n+1\}.$$
Hence, we have
$$ [BY_2]=\sum_{k=0}^{n+1} \left(\alpha_kI_m+\beta_k\exp (-A(b-a))\right)(-A)^k.$$
Therefore, the characteristic block matrix becomes
$$ M(L, B)= \left( \alpha_0+\beta_0; \sum_{k=0}^{n+1} \left(\alpha_k+\beta_k\exp (-A(b-a))\right)(-A)^k\right).$$

According to Theorem \ref{th dimker}, the dimensions of the kernel and co-kernel of the inhomogeneous boundary-value problem are equal, resp., to the dimensions of the kernel and co-kernel of the matrix $M(L, B)$.

In particular, if $\beta_k \equiv 0$ and the problem is one-point, then the block characteristic matrix takes the form
$$ M(L, B)= \left( \alpha_0; \sum_{k=0}^{n+1} \alpha_k(-A)^k\right).$$

Therefore, in this case, the $d$-characteristics of the boundary-value problem do not depend on the length of the interval $(a,b)$.

Note that the matrix $\exp (-A(b-a))$ can be found in an explicit form since every entire analytic function of the number matrix $A \in \mathbb{C}^{m \times m}$ is a polynomial of $A$. This polynomial is expressed via the matrix $A$ by the Lagrange--Sylvester Interpolation Formula (see, e.g., \cite{Gan1959}). Its degree is no greater than $m-1$.

\textbf{Example 4.} Consider a two-point boundary-value problem for another system of second-order differential equations
\begin{equation*}
    (Ly)(t):= y^{\prime \prime} (t)+Ay(t),\quad
t \in(a,b),
\end{equation*}
where $A\in\mathbb{C}^{m\times m}$. The boundary conditions induced by the same operator as that in Example~3; namely
\begin{equation*}
By=\sum^{n+1} _{k=0} \left(\alpha_{k} y^{(k)}(a)+ \beta_{k}y^{(k)}(b)\right).
\end{equation*}

It is easy to check in this case that, for each fixed $t \in [a, b]$, the fundamental matrix-valued functions $Y_{1}(t)$ and $Y_{2}(t)$ are entire functions of the number matrix $A$ given by some convergent power series. Then
\begin{gather*}
[BY_1]=\sum_{\substack{k=1\\k \,\, is \,\, odd}}^{n+1} \beta_k(-1)^kA^{k}
+\sum_{\substack{k=0\\k \,\, is \,\, even}}^{n+1} \alpha_k(-1)^k(\sqrt{A})^{2k-1}\sin \big(\sqrt{A}(b-a)\big)\\
+\sum_{\substack{k=1\\k \,\, is \,\, odd}} \beta_k(-1)^kA^{k}\cos \big(\sqrt{A}(b-a)\big)
\end{gather*}
and
\begin{gather*}
[BY_2]=\sum_{\substack{k=1\\k \,\, is \,\, even}}^{n+1} \alpha_k(-1)^kA^{k}
+\sum_{\substack{k=0\\k \,\, is \,\, even}}^{n+1} \alpha_k(-1)^kA^{k}\cos \big(\sqrt{A}(b-a)\big)\\
+\sum_{\substack{k=1\\k \,\, is \,\, odd}}^{n+1} \beta_k(-1)^k(\sqrt{A})^{2k-1}\sin \big( \sqrt{A}(b-a)\big),
\end{gather*}
with the block characteristic matrix $M(L,B)=[BY_1; BY_2]$.

Specifically, if $\beta_k \equiv 0$ (the case of the one-point boundary-value problem), then
\begin{gather*}
M(L,B) = \left[\,\sum_{\substack{k=0\\k \,\, is \,\, even}}^{n+1} \alpha_k(-1)^k(\sqrt{A})^{2k-1}\sin \big(\sqrt{A}(b-a)\big); \right. \\
\left.\sum_{\substack{k=1\\k \,\, is \,\, even}}^{n+1} \alpha_k(-1)^kA^{k}+
\sum_{\substack{k=0\\k \,\, is \,\, even}}^{n+1}
\alpha_k(-1)^kA^{k}\cos \big(\sqrt{A}(b-a)\big)  \right].
\end{gather*}
Unlike Example 3, this matrix depends in general on the length of the interval $(a,b)$.

If $\alpha_k \equiv 0$, $k$ is even, $\beta_k \equiv 0$, and $k$ is odd, then the characteristic matrix $M(L,B)=O_{2m \times l}$. Therefore, its Fredholm numbers take the largest possible values.

As in Example 3, the matrices $\sin \big(\sqrt{A}(b-a)\big)$ and $\cos \big( \sqrt{A}(b-a)\big)$ can be exactly found as Lagrange--Sylvester interpolation polynomials.

\textbf{Example 5.} Consider the following linear boundary-value problem for a system of $m$ first-order differential equations:
\begin{equation}\label{1.6.1t}
    Ly(t):= y'(t)=f(t),\quad
t \in(a,b), \quad By= c,
\end{equation}
where $f(\cdot) \in (W_{p}^{n})^{m}$, and $c \in \mathbb{C}^{l}$, and $B$
is an arbitrary linear continuous operator from
$(W_{p}^{n+1})^{m}$ to $\mathbb{C}^{l}$. We suppose that $1\leq p<\infty$.

Note that $Y(\cdot)=I_{m}$ is the unique solution of the linear homogeneous matrix equation of the form \eqref{1.6.1t} with the initial Cauchy condition
\begin{equation*}\label{r3}
   Y'(t)=0,\quad t\in (a,b), \qquad Y(a)=I_{m}.
  \end{equation*}
According to \eqref{st anal}, we have
\begin{equation*}
M(L,B)=[BY]=\sum _{i=0}^{n} \alpha_{i} Y^{(i)}(a)+
\int_{a}^b \Phi(t)Y^{(n+1)}(t){\rm d}t=\alpha_{0}.
\end{equation*}
Therefore,
\begin{equation*}
\operatorname{dim} \operatorname{ker}(M(L,B))=\operatorname{dim} \operatorname{ker}(\alpha_{0})
\end{equation*}
\begin{equation*}
\operatorname{dim} \operatorname{coker}(M(L,B))=\operatorname{dim} \operatorname{coker}(\alpha_{0}).
\end{equation*}
Hence, the boundary-value problem \eqref{1.6.1t} is well posed if and only if the number matrix $\alpha_{0}$ is square and nonsingular.

The results given in this section were obtained in \cite{AtlMikh2024Zh, AtlMikh2024}. In Sobolev-Slobodetskii spaces, Example~1 and a special case of Example~2 (two-point problem) are given in \cite{MikhAtlSkor}.

\section{Convergence of the characteristic matrices}\label{section4}

Together with the problem \eqref{bound_pr_1}, \eqref{bound_pr_2}, we consider the sequence of boun\-da\-ry-value problems
\begin{gather}\label{6.syste}
L(k)y(t,k):= y^{(r)}(t,k)+\sum_{j=1}^{r}A_{r-j}(t,k)y^{(r-j)}(t,k)=f(t,k),\\
B(k)y(\cdot,k)=c(k), \quad t\in (a,b), \quad k\in\mathbb{N},\label{6.kue}
\end{gather}
where the matrix-valued functions $A_{r-j}(\cdot,k)$, the vector-valued functions $f(\cdot,k)$, the vectors $c(k)$, and the linear continuous operators $B(k)$ satisfy the above conditions imposed on the problem \eqref{bound_pr_1}, \eqref{bound_pr_2}. We assume in the sequel that $k \in \mathbb{N}$ and that all asymptotic relations are considered for $k\rightarrow \infty$. The boundary-value problem \eqref{6.syste}, \eqref{6.kue} is also the most general (generic) with respect to the Sobolev space $W^{n+r}_p$.

We associate the sequence of linear continuous operators
\begin{equation}\label{L(e)_B(e)}
(L(k),B(k))\colon(W^{n+r}_p)^{m}\rightarrow (W^{n}_p)^{m}\times\mathbb{C}^{l}
\end{equation}
and the sequence of characteristic matrices
$$
M\big(L(k),B(k)\big):=
\big(\left[B(k)Y_1(k)\right],\dots,\left[B(k)Y_{r}(k)\right]\big) \subset \mathbb{C}^{l\times rm}
$$
with the boundary-value problems \eqref{6.syste}, \eqref{6.kue}.

As usual,
\begin{equation}\label{zb LB11}
\left(L(k),B(k)\right)\xrightarrow{s} \left(L,B\right)
\end{equation}
denotes the strong convergence of the sequence of operators $(L(k),B(k))$ to the operator~$(L,B)$.

The next theorem provides a sufficient condition for the convergence of the sequence of characteristic matrices $M\left(L(k),B(k)\right)$ to the matrix $M\left(L,B\right)$.

\begin{theorem}\label{koef matr}
If the sequence of operators $(L(k),B(k))$ converges strongly to the operator $(L,B)$, then the sequence of characteristic matrices $M\big(L(k),B(k)\big)$ converges to the matrix $M\big(L,B\big)$; i.e.,
\begin{equation*}\label{zb har m}
(L(k),B(k))\xrightarrow{s} (L,B)\;\Longrightarrow\; M(L(k),B(k))\rightarrow M(L,B).
\end{equation*}
\end{theorem}

Theorem \ref{koef matr} implies

\begin{theorem}\label{ker coker}
If condition \eqref{zb LB11} is satisfied, then the following inequalities hold true for all sufficiently large $k$:
\begin{gather*}
\operatorname{dim} \operatorname{ker}(L(k),B(k))\leq \operatorname{dim} \operatorname{ker}(L,B), \label{ner ker} \\
\operatorname{dim} \operatorname{coker}(L(k),B(k))\leq \operatorname{dim} \operatorname{coker}(L,B). \label{ner coker}
\end{gather*}
\end{theorem}

Let us consider three significant direct consequences of Theorem~\ref{ker coker}. Suppose that condition \eqref{zb LB11} is satisfied.

\begin{corollary}\label{cor1}
If the operator $(L,B)$ is invertible, then so are the operators $\left(L(k),B(k)\right)$ for all sufficiently large $k$.
\end{corollary}

\begin{corollary}\label{cor2}
If the boundary-value problem \eqref{bound_pr_1}, \eqref{bound_pr_2} has a  solution for arbitrarily chosen right-hand sides, then so do the boundary-value problems \eqref{6.syste}, \eqref{6.kue} for all sufficiently large $k$.
\end{corollary}

\begin{corollary}\label{cor3}
If the homogeneous boundary-value problem \eqref{bound_pr_1}, \eqref{bound_pr_2} has only a trivial solution, then so do the homogeneous problems \eqref{6.syste}, \eqref{6.kue} for all sufficiently large $k$.
\end{corollary}

Note that the conclusion of Theorem \ref{ker coker} and its consequences cease to be valid for arbitrary bounded linear operators between infinite-dimensional Banach spaces.

The results presented in this section were obtained in \cite{AtlMikh2024}.

\section{Continuity of solutions in a parameter}\label{section5}

Let us consider the linear boundary-value problem
\begin{equation}\label{eq7}
L(\varepsilon)y(t,\varepsilon):=y^{(r)}(t,\varepsilon) + \sum\limits_{j=1}^rA_{r-j}(t,\varepsilon)y^ {(r-j)}(t,\varepsilon)=f(t,\varepsilon),
\end{equation}
\begin{equation}\label{eq8}
 B(\varepsilon)y(\cdot;\varepsilon) = c(\varepsilon), \quad t\in (a,b),
\end{equation}
parameterized by number $\varepsilon \in [0,\varepsilon_0)$, $\varepsilon_0>0$. Here, for every fixed $\varepsilon$, the matrix-valued functions $A_{r-j}(\cdot;\varepsilon) \in (W^{n}_p)^{m\times m}$,  vector-valued function $f(\cdot;\varepsilon) \in (W^{n}_p)^m$, vector $c(\varepsilon)\in\mathbb{C}^{rm}$, and the linear continuous operator $$
B(\varepsilon) \colon (W^{n+r}_p)^m\rightarrow\mathbb{C}^{rm}
$$
are given, whereas the vector-valued function $y(\cdot;\varepsilon) \in (W^{n+r}_p)^m$ is unknown.

It follows from Theorem \ref{th_fredh high} that the boundary-value problem \eqref{eq7}, \eqref{eq8} is a Fredholm one with index zero.

\begin{definition}\label{d2}
The solution to the boundary-value problem \eqref{eq7}, \eqref{eq8} depends continuously on the parameter $\varepsilon$ at $\varepsilon=0$ if the following two conditions are satisfied:
\begin{itemize}
\item [$(\ast)$] there exists a positive number $\varepsilon_{1}<\varepsilon_{0}$ such that, for any $\varepsilon\in[0,\varepsilon_{1})$ and an arbitrary chosen right-hand sides $f(\cdot;\varepsilon)\in (W^{n}_p)^{m}$ and $c(\varepsilon)\in\mathbb{C}^{rm}$, this problem has a unique solution $y(\cdot;\varepsilon)$ that belongs to the space $(W^{n+r}_p)^{m}$;

\item [$(\ast\ast)$] the convergence of the right-hand sides $f(\cdot;\varepsilon)\to f(\cdot;0)$ in $(W_p^{n})^{m}$ and $c(\varepsilon)\to c(0)$ in $\mathbb{C}^{rm}$ implies the convergence of the solutions $y(\cdot;\varepsilon)\to y(\cdot;0)$ in $(W^{n+r}_p)^{m}$.
\end{itemize}
\end{definition}

Here and below, the limits are considered as $\varepsilon\to0+$.

Definition \ref{d2} is equivalent to the following two conditions:
\begin{itemize}
\item [---] The operator $\big(L(\varepsilon), B(\varepsilon)\big)$ is invertible for all sufficiently small $\varepsilon>0$;
\item [---] $\big(L(\varepsilon), B(\varepsilon)\big)^{-1} \stackrel{s}{\longrightarrow} \big(L(0), B(0)\big)^ {-1}$.
\end{itemize}

Consider the following assumptions:
\begin{itemize}
  \item [(0)] the homogeneous boundary-value problem has only the trivial solution
$$L(0)y(t,0)=0,\quad t\in(a, b),\quad B(0)y(\cdot,0)=0;$$
  \item [(I)] $A_{r-j}(\cdot;\varepsilon)\to A_{r-j}(\cdot;0)$ in the space $(W^{n}_p)^{m\times m}$ for each number $j\in\{1,\ldots, r\}$;
  \item [(II)] $B(\varepsilon)y\to B(0)y$ in the space $\mathbb{C}^{rm}$ for every $y\in(W^{n+r}_p)^m$.
\end{itemize}

\begin{theorem}\label{th4} The solution to the boundary-value problem \eqref{eq7}, \eqref{eq8} depends continuously on the parameter $\varepsilon$ at $\varepsilon=0$ if and only if this problem satisfies conditions $(0)$, $(\mathrm{I})$, and $(\mathrm{II})$.
\end{theorem}

This Theorem implies that, if the operator $\big(L(0), B(0)\big)$ is invertible, then
$$\big(L(\varepsilon), B(\varepsilon)\big) \stackrel{s}{\longrightarrow} \big(L(0), B(0)\big) \Leftrightarrow \big(L(\varepsilon), B(\varepsilon)\big)^{-1} \stackrel{s}{\longrightarrow} \big(L(0), B(0)\big)^{-1}.$$
Note that the conclusion of Theorem \ref{th4} and its consequences cease to be valid for arbitrary bounded linear operators between infinite-dimensional Banach spaces. Note that the set of all irreversible operators is everywhere dense in the strong operator topology.

We supplement our result with a two-sided estimate of the error
$\bigl\|y(\cdot;0)-y(\cdot;\varepsilon)\bigr\|_{n+r,p}$ of the solution $y(\cdot;\varepsilon)$ via its discrepancy
$$
\widetilde{d}_{n,p}(\varepsilon):=
\bigl\|L(\varepsilon)y(\cdot;0)-f(\cdot;\varepsilon)\bigr\|_{n,p}+
\bigl\|B(\varepsilon)y(\cdot;0)-c(\varepsilon)\bigr\|_{\mathbb{C}^{rm}}.
$$
Here, we interpret $y(\cdot;\varepsilon)$ as an approximate solution to the problem \eqref{eq7}, \eqref{eq8} with $\varepsilon = 0$.

\begin{theorem}\label{th5} Suppose that the boundary-value problem \eqref{eq7}, \eqref{eq8} satisfies conditions $(0)$, $(\mathrm{I})$, and $(\mathrm{II})$. Then there exist positive numbers $\varepsilon_{2}<\varepsilon_{1}$ and $\gamma_{1}$, $\gamma_{2}$ such that, for any  $\varepsilon\in(0,\varepsilon_{2})$, the following two-sided estimate is true:
$$
\gamma_{1}\,\widetilde{d}_{n,p}(\varepsilon)
\leq\bigl\|y(\cdot;0)-y(\cdot;\varepsilon)\bigr\|_{n+r,p}\leq
\gamma_{2}\,\widetilde{d}_{n,p}(\varepsilon),
$$
where the numbers $\varepsilon_{2}$, $\gamma_{1}$, and $\gamma_{2}$ do not depend of $y(\cdot;\varepsilon)$ and $y(\cdot;0)$.
\end{theorem}

Thus, the error and discrepancy of the solution $y(\cdot;\varepsilon)$ to the boundary-value problem \eqref{eq7}, \eqref{eq8} are of the same degree of smallness.

The results presented in this section were obtained in \cite{AtlMikh2024Zh}. Unlike the method used in \cite{GKM2017}, our approach is more general and allows studying the solutions of boundary-value problems not only in Sobolev spaces, but also in other function spaces (see, e.g., \cite{MMS2016}). The case of first order equations (with $r=1$) in Sobolev spaces is considered in \cite{Atl2}, and in the case of Sobolev-Slobodetskii spaces in \cite{Hnyp}. For the most general class of inhomogeneous boundary-value problems for systems of differential equations of an arbitrary order whose solutions belong to the Sobolev space with $1\leq p < \infty$ similar results were obtained in \cite{GnypKodMik15}; the case where the solutions range over an appropriate Sobolev--Slobodetskii space was studied in \cite{MaslyukMykhailets2018}.

\section{Multipoint problems} \label{section6}
The results of this section are the principle of averaging for solutions of problems with multipoint boundary conditions.

We consider the most general class of multipoint linear boundary-value problems for systems of ordinary differential equations of any order whose solutions belong to the Sobolev space $W_p^{n+r}$. We consider the case where the points of the closed interval $[a,b]$ appearing in boundary conditions are not fixed and depend on a number parameter and the number of these points may change. The case
$p=\infty$ is special and has not been studied earlier.

We arbitrarily choose $N$ different points $\{t_1,\ldots,t_{N}\} \subset [a, b]$ and consider a multipoint boundary-value problem of the form
\begin{gather}\label{bound_pr_1-bis}
Ly(t)\equiv y^{(r)}(t)+ \sum\limits_{j=1}^{r}A_{r-j}(t)y^{(r-j)}(t)=
f(t),\quad t \in(a,b),\\
By\equiv\sum_{l=0}^{n+r-1}\sum_{j=1}^{N}
\beta_{j}^{(l)}y^{(l)}(t_j)=q, \label{bagat_cond_}
\end{gather}
where $y\in(W^{n+r}_p)^m$ is an unknown vector function, $A_{r-j}\in(W_p^n)^{m\times m}$ are arbitrary matrix functions, $f\in(W^n_p)^m$ are arbitrary vector functions, $\beta_{j}^{(l)}\in\mathbb{C}^{rm\times m}$ are arbitrary matrices, and $q\in\mathbb{C}^{rm}$.

In view of the continuous embedding
\begin{equation}\label{W-C-embedding}
(W^{n+r}_{p})^m\hookrightarrow(C^{n+r-1})^m,
\end{equation}
the left-hand side of the boundary condition \eqref{bagat_cond_} makes sense, and the mapping $y\mapsto By$, where $y\in(W^{n+r}_p)^m$, is a continuous operator from the space $(W^{n+r}_p)^m$ to $\mathbb{C}^{rm}$. Note that the boundary condition \eqref{bagat_cond_} is not classical because it contains the derivatives $y^{(l)}$ of integral orders $l$, where $0 < l \leq n +r-1$.

Problem \eqref{bound_pr_1-bis}, \eqref{bagat_cond_} is regarded as the limit boundary-value problem, as $\varepsilon\to0+$, for the following multipoint
boundary-values problem depending on the parameter $\varepsilon\in(0,\varepsilon_0)$:
\begin{gather}
L(\varepsilon)y(t,\varepsilon):=  y^{(r)}(t,\varepsilon) + \sum\limits_{j=1}^rA_{r-j}(t,\varepsilon)y^{(r-j)}(t,\varepsilon) =f(t,\varepsilon),\quad t \in(a,b),\label{7syste} \\
B(\varepsilon)y(\cdot,\varepsilon)= \sum\limits_{j=0}^{N}\sum\limits_{k=1}^{\omega_j(\varepsilon)}
\sum\limits_{l=0}^{n+r-1}{\beta_{j,k}^{(l)}
(\varepsilon)y^{(l)}(t_{j,k}(\varepsilon),\varepsilon)}=q(\varepsilon). \label{7kue}
\end{gather}
Here, for any fixed value of the parameter $\varepsilon$, the vector function  $y(\cdot,\varepsilon)\in(W^{n+r}_p)^m$ is unknown and the matrix functions $A_{r-j}(\cdot,\varepsilon)\in(W_p^n)^{m\times m}$, vector function $f(\cdot,\varepsilon)\in(W^n_p)^m$, vectors $q(\varepsilon)\in \mathbb{C}^{rm}$, and matrices $\beta_{j,k}^{(l)}(\varepsilon)\in \mathbb{C}^{m\times m}$ are given. For $\varepsilon>0$, we take at least $N$ points $t_{j,k}(\varepsilon)$ of the segment gathered into $N+1$ series as follows: for any fixed $j\in\{1,\ldots,N\}$, all points $t_{j,k}(\varepsilon)$ must have the same limit $t_j$ as $\varepsilon \rightarrow 0+$. This requirement is not imposed on the
points $t_{0,k}(\varepsilon)$. Note that the series with $j=0$ may be absent.

We represent vectors and vector-valued functions in the form of columns. The solution of the boundary-value problem \eqref{7syste}, \eqref{7kue} is defined to be a vector-valued function $y(\cdot,\varepsilon)\in (W_{p}^{n+r})^m$ satisfying both equation $(\ref{7syste})$ everywhere (for $n\geq 1$) and almost everywhere (for $n= 0$) on $(a,b)$ and equality~\eqref{7kue} specifying $rm$ scalar boundary conditions. The presence of the repeated sum over the indices $j$ and $k$ in the boundary condition \eqref{7kue} is explained by the subsequent assumptions concerning the behavior of the points $t_{j,k}(\varepsilon)$ as $\nobreak{\varepsilon\to0+}$ depending on~$j$.

In the limit case of $\varepsilon=0$, we consider a boundary-value problem
\begin{gather}
\label{7lse}L(0)y(t,0)=f(t,0),\quad t \in(a,b),\\
\label{7lkue} B(0)y(\cdot,0)=\sum\limits_{j=1}^{N}
\sum\limits_{l=0}^{n+r-1}{\beta_{j}^{(l)}y^{(l)}(t_{j},0)}=q(0),
\end{gather}
where the matrices $\beta_{j}^{(l)}\in\mathbb{C}^{m\times m}$, the points $t_{j}\in[a,b]$, and the vector $q(0)\in \mathbb{C}^{rm}$ are arbitrary.

For any $\varepsilon\in[0,\varepsilon_0)$, $B(\varepsilon)$ is a continuous linear operator
\begin{equation}\label{7B}
B(\varepsilon)\colon (W_{p}^{n+r})^m \to \mathbb{C}^{rm}.
\end{equation}

For every $\varepsilon\in[0,\varepsilon_0)$, the boundary-value problem \eqref{7syste}, \eqref{7kue} induces the linear operator
\begin{equation}\label{operatop3.77}
\big(L(\varepsilon),B(\varepsilon)\big)\colon (W^{n+r}_p)^m \rightarrow (W^{n}_p)^m\times\mathbb{C}^{rm}.
\end{equation}
According to Theorem \ref{th_fredh high}, \eqref{operatop3.77}  is a bounded Fredholm operator with zero index.

Let us establish explicit sufficient conditions for the solution $y=y(\cdot,\varepsilon)$ of the multipoint boundary-value problem \eqref{7syste}, \eqref{7kue} to be continuous with respect to the parameter $\varepsilon$ in the Sobolev space $W^{n+r}_p$, with $1\leq p\leq \infty$; i.e., the conditions for the solution $y(\cdot,\varepsilon)$ to exist, be unique, and satisfy the limit relation
\begin{equation}\label{7zb r}
\bigl\|y(\cdot,\varepsilon)-y(\cdot,0) \bigr\|_{n+r,p} \to 0 \quad \mbox{as} \quad \varepsilon\to0+.
\end{equation}

In order that this task make sense, we assume the following:

\textbf{Condition (0).} \emph{The homogeneous limit boundary-value problem of the form \eqref{7lse}, \eqref{7lkue} has only the trivial solution, i.e., is not degenerate.}

This implies that, for $\varepsilon = 0$, the Fredholm operator \eqref{operatop3.77} is an isomorphism, i.e.
$$
\big(L(0),B(0)\big)\colon (W^{n+r}_p)^m \leftrightarrow (W^{n}_p)^m\times\mathbb{C}^{rm}.
$$

Hence, the boundary-value problem \eqref{7lse}, \eqref{7lkue} has a unique solution $y(t,0)\in(W^{n+r}_p)^m$ for any right-hand sides $f(t,0)\in(W^{n}_p)^m$ and $q(0) \in \mathbb{C}^{rm}$.

We consider the following

\textbf{Assumptions} as $\varepsilon\to0+$:
\begin{itemize}
    \item [$(\alpha)$] $t_{j,k}(\varepsilon)\to t_{j}$ for all  $j\in\{1,\ldots,N\}$ and $k\in\{1,\ldots,\omega_j(\varepsilon)\};$

    \item [$(\beta)$] $\sum\limits_{k=1}^{\omega_j(\varepsilon)}
\beta_{j,k}^{(l)}(\varepsilon)\to\beta_{j}^{(l)}$ for all $j\in\{1,\ldots,N\}$ and $l\in\{0,\ldots,n+r-1\}$;

    \item [$(\gamma)$] $\sum\limits_{k=1}^{\omega_j(\varepsilon)}\left\|\beta_{j,k}^{(l)}(\varepsilon)\right\|
\left|t_{j,k}(\varepsilon)-t_j\right|\to0$ for all  $j\in\{1,\ldots,N\}$, $k\in\{1,\ldots,\omega_j(\varepsilon)\}$, and $l\in\{0,\ldots,n+r-1\}$;

      \item [$(\delta)$] $\sum\limits_{k=1}^{\omega_0(\varepsilon)}\left\|\beta_{0,k}^{(l)}(\varepsilon)\right\|\to0$  for all $k\in\{1,\ldots,\omega_0(\varepsilon)\}$ and $l\in\{0,\ldots,n+r-1\}$.
\end{itemize}

Note that, for the boundary-value problem \eqref{7syste}, \eqref{7kue}, we do not suppose that the coefficients $A_{r-j}(\cdot,\varepsilon)$ and $\beta_{j,k}^{(l)}(\varepsilon)$ or the points $t_{j,k}(\varepsilon)$ are characterized by a certain regularity with respect to the parameter $\varepsilon >0$. Assume that, for any fixed $j\in\{1,\ldots,N\}$, all points   $t_{j,k}(\varepsilon)$ have the same limit as $\varepsilon\to0+$. At the same time, this requirement is not imposed on the points of the zero series $t_{0,k}(\varepsilon)$.

In the conditions $(\gamma)$ and $(\delta)$, the expression $\|\cdot\|$ denotes a norm of a complex number matrix; this norm is equal to the sum of absolute values of all elements of the matrix. In view of assumptions $(\beta)$ and $(\gamma)$, the norms of the coefficients $\big\|\beta_{j,k}^{(l)}(\varepsilon)\big\|$ may infinitely (but not very rapidly) increase as $\varepsilon\to0+$. It follows from
condition $(\delta)$ that, unlike condition $(\alpha)$, it is not necessary to demand the convergence of the points $t_{0,j}(\varepsilon)$ as $\varepsilon\to0+$.

Let us formulate limit theorems for the solutions to the multipoint boundary-value problem \eqref{7syste}, \eqref{7kue} in the $p=\infty$
case.

\begin{theorem}\label{6.dost_t} Suppose that the boundary-value problem \eqref{7syste}, \eqref{7kue} with $p=\infty$ satisfies the assumptions $(\alpha)$, $(\beta)$, $(\gamma)$, and $(\delta)$. Then it satisfies the limit condition~\textup{(II)}. Moreover, if conditions \textup{(0)} and \textup{(I)} are satisfied, then, for sufficiently small $\varepsilon$, the solution to this problem exists, is unique, and satisfies the limit relation \eqref{7zb r}.
\end{theorem}

Focusing now on the case $1\leq p< \infty$, we consider the following

\textbf{Assumptions} as $\varepsilon\to0+$:
\begin{itemize}
       \item[$(\gamma_p)$] $\sum\limits_{k=1}^{\omega_j(\varepsilon)}\left\|\beta_{j,k}^{(n+r-1)}(\varepsilon)\right\|
\left|t_{j,k}(\varepsilon)-t_j\right|^{1/p'}=O(1)$ for all $j\in\{1,\ldots,N\}$ and $k\in\{1,\ldots,\omega_j(\varepsilon)\}$, where $1/p+1/p'=1$;

    \item[$(\gamma')$] $\sum\limits_{k=1}^{\omega_j(\varepsilon)}\left\|\beta_{j,k}^{(l)}(\varepsilon)\right\|
\left|t_{j,k}(\varepsilon)-t_j\right|\to0$ for all $j\in\{1,\ldots,N\}$, $k\in\{1,\ldots,\omega_j(\varepsilon)\}$, and
 $l\in\{0,\ldots,n+r-2\}$.
    \end{itemize}

Note that the systems of conditions $(\alpha)$, $(\beta)$, $(\gamma)$, $(\delta)$ and $(\alpha)$, $(\beta)$, $(\gamma_p)$, $(\gamma')$, $(\delta)$ do not guarantee the uniform convergence of continuous operators $B(\varepsilon):(W^{n+r}_p)^m\to \mathbb{C}^{rm}$ to $B(0)$ as $\varepsilon\to0+$. For this reason, Theorem~\ref{6.dost_t} does not follow from the general facts of the theory of linear operators.

\begin{theorem}\label{th4.1}
Let $1\leq p< \infty$, and suppose that the boundary-value problem \eqref{7syste}, \eqref{7kue} satisfies assumptions $(\alpha)$, $(\beta)$, $(\gamma_p)$, $(\gamma')$, and $(\delta)$. Then this problem satisfies the limit condition~\textup{(II)}. Moreover, if conditions \textup{(0)} and \textup{(I)}  are satisfied, then, for sufficiently small $\varepsilon$, the solution of the problem exists, is unique, and satisfies the limit relation \eqref{7zb r}.
\end{theorem}

Papers \cite{Atl3,Atl4} gives sufficient conditions for the continuity of solutions to multipoint boundary-value problems with respect to the parameter in Sobolev spaces.

Note that, for first order differential equations ($r=1$), Theorems \ref{6.dost_t} and \ref{th4.1} were proved in \cite{Atl3}. In the general case, for differential equations of any order, the proofs of Theorems \ref{6.dost_t} and \ref{th4.1} are based on the criterion of continuity of the most general boundary-value problems \cite{AtlMikh2019Dop}.

Note that papers \cite{Kodlyuk2012Coll, Kodliuk2012Dop}  investigated multipoint boundary-value problems for systems of first order ordinary differential equations in the Sobolev spaces $W^{n}_p$, where $1\leq p< \infty$. However, in these papers, the points of the segment $[a,b]$ appearing in the limit condition are fixed and independent of the parameter. Paper \cite{Hnyp-Kodl15} studied nonclassical multipoint boundary-value problems for systems of ordinary differential equations of arbitrary order in the Sobolev spaces  $W^{n+r}_p$, where $1\leq p< \infty$. However, in this paper, the number of points in each series is independent of the parameter~$\varepsilon$.

\section{Acknowledgments}

\textit{The authors thank Prof. Aleksandr Murach for his discussion of the paper and valuable remarks.}

The work of the first named author was funded by the National Academy of Sciences of Ukraine and by a grant from the Simons Foundation (1290607, V.A.M.). The work of the second named author was funded by the Academy of Finland within grant no.~359642 by the Research Council of Finland: Researcher mobility to Finland. The authors was also supported by the European Union's Horizon 2020 research and innovation programme under the Marie Sk{\l}odowska-Curie grant agreement No~873071 (SOMPATY: Spectral Optimization: From Mathematics to Physics and Advanced Technology).


\textit{Mikhailets V.A.:} Institute of Mathematics of the National Academy of Sciences of Ukraine, st. Tereschenkivska 3, 01024 Kyiv, Ukraine \\ mikhailets@imath.kiev.ua\\

\textit{Atlasiuk O.M.} Institute of Mathematics of the National Academy of Sciences of Ukraine, st. Tereschenkivska 3, 01024 Kyiv, Ukraine and University of Helsinki, st.  Yliopistonkatu 4, 00100 Helsinki, Finland \\ hatlasiuk@gmail.com

\end{document}